\documentclass[12pt]{article}
\usepackage{amsfonts}
\usepackage{amssymb}
\usepackage{amsmath}
\usepackage{hyperref}
\usepackage{color}
=0pt
\def\sqr#1#2{{\vcenter{\vbox{\hrule height.#2pt
              \hbox{\vrule width.#2pt height#1pt \kern#1pt \vrule
width.#2pt}
              \hrule height.#2pt}}}}
\def\signed #1{{\unskip\nobreak\hfil\penalty50
              \hskip2em\hbox{}\nobreak\hfil#1
              \parfillskip=0pt \finalhyphendemerits=0 \par}}
\def\endpf{\signed {$\sqr69$}}
\def\dbE{{\mathbb{E}}}
\def\dbF{{\mathbb{F}}}

\def\dbN{{\mathbb{N}}}

\def\dbP{{\mathbb{P}}}

\def\dbR{{\mathbb{R}}}

\def\a{\alpha}
\def\b{\beta}

\def\d{\delta}
\def\e{\varepsilon}

\def\t{\tau}
\def\f{\varphi}
\def\th{\theta}

\def\3n{\negthinspace \negthinspace \negthinspace }
\def\2n{\negthinspace \negthinspace }
\def\1n{\negthinspace }
\def\ns{\noalign{\smallskip} }

\def\ds{\displaystyle}
\def\Th{\Theta}
\def\L{\Lambda}

\def\Om{\Omega}
\def\om{\omega}
\def\cA{{\cal A}}
\def\cB{{\cal B}}

\def\cF{{\cal F}}

\def\cJ{{\cal J}}

\def\cR{{\cal R}}
\def\cS{{\cal S}}

\def\mE{{\mathbb{E}}}

\def\no{\noindent}

\def\ms{\medskip}
\def\bs{\bigskip}
\def\q{\quad}
\def\qq{\qquad}
\def\hb{\hbox}

\def\wt{\widetilde}
\def\cd{\cdot}

\def\ae{\hbox{\rm a.e.{ }}}
\def\as{\hbox{\rm a.s.{ }}}

\def\({\Big (}
\def\){\Big )}
\def\[{\Big[}
\def\]{\Big]}

\def\={\buildrel \triangle \over =}

\def\be{\begin{equation}}
\def\bel{\begin{equation}\label}
\def\ee{\end{equation}}
\def\bea{\begin{eqnarray}}
\def\eea{\end{eqnarray}}
\def\bt{\begin{theorem}}
\def\et{\end{theorem}}
\def\bc{\begin{corollary}}
\def\ec{\end{corollary}}
\def\bl{\begin{lemma}}
\def\el{\end{lemma}}
\def\bp{\begin{proposition}}
\def\ep{\end{proposition}}
\def\br{\begin{remark}}
\def\er{\end{remark}}
\def\ba{\begin{array}}
\def\ea{\end{array}}
\def\bd{\begin{definition}}
\def\ed{\end{definition}}

\newtheorem{lemma}{Lemma}[section]
\newtheorem{remark}{Remark}[section]
\newtheorem{example}{Example}[section]
\newtheorem{theorem}{Theorem}[section]
\newtheorem{corollary}{Corollary}[section]

\newtheorem{definition}{Definition}[section]
\newtheorem{proposition}{Proposition}[section]

\makeatletter
   
   \@addtoreset{equation}{section}
\makeatother

\begin{document}

\title{\bf Characterization of Optimal Feedback for Stochastic Linear Quadratic Control Problems}

\author{Qi L\"u\thanks{School of
Mathematics, Sichuan University, Chengdu
610064, Sichuan Province, China.  {\small\it E-mail:} {\small\tt
lu@scu.edu.cn}.\ms}, \q Tianxiao Wang\thanks{School of
Mathematics, Sichuan University, Chengdu
610064, Sichuan Province, China.  {\small\it E-mail:} {\small\tt
xiaotian2008001@gmail.com}.\ms} \q and \q Xu Zhang\thanks{
School of Mathematics, Sichuan University,
Chengdu 610064, Sichuan Province, China.  {\small\it E-mail:}
{\small\tt zhang\_xu@scu.edu.cn}.} }

\date{}

\maketitle

\begin{abstract}
One of the fundamental issues in Control Theory
is to design feedback controls. It is well-known
that, the purpose of introducing Riccati
equations in the deterministic case is to
provide the desired feedback controls for linear
quadratic control problems. To date, the same
problem in the stochastic setting is only
partially well-understood. In this paper, we
establish the equivalence between the existence
of optimal feedback controls for the stochastic
linear quadratic control problems with random
coefficients and the solvability of the
corresponding backward stochastic Riccati
equations  in a suitable sense. We also give a counterexample showing the nonexistence of feedback controls to a solvable stochastic
linear quadratic control problem. This is a new phenomenon in the stochastic setting, significantly different from its deterministic counterpart.

\end{abstract}

\bs

\no{\bf 2010 Mathematics Subject
Classification}. Primary  93E20; Secondary  93B52, 93C05, 60H10.

\bs

\no{\bf Key Words}. Stochastic linear quadratic
problem, feedback control, backward stochastic
Riccati equation,  backward stochastic
differential equation.

\section{Introduction}\label{s1}

Let $T>0$ and $(\Om,\cF,\dbF,\dbP)$ be a complete filtered
probability space with
$\dbF=\{\cF_t\}_{t\in[0,T]}$, which is the natural
filtration generated by a one-dimensional  standard
Brownian motion $\{W(t)\}_{t\in[0,T]}$.

For any
$k\in\dbN$, $t\in[0,T]$ and $r\in [1,\infty)$,
denote by $L_{\cF_t}^r(\Om;\dbR^k)$ the Banach
space of all $\cF_t$-measurable random variables
$\xi:\Om\to \dbR^k$ so that
$\mathbb{E}|\xi|_{\dbR^k}^r < \infty$, with the
canonical  norm. Denote by
$L^{r}_{\dbF}(\Om;C([t,T];\dbR^k))$ the Banach
space of all $\dbR^k$-valued $\dbF$-adapted,
continuous stochastic processes $\phi(\cdot)$,
with the following norm
$$
|\phi(\cd)|_{L^{r}_{\dbF}(\Om;C([t,T];\dbR^k))}
\= \(\mE\sup_{\tau\in
[t,T]}|\phi(\tau)|_{\dbR^k}^r\)^{1/r}.
$$
Fix any $r_1,r_2,r_3,r_4\in[1,\infty)$. Put
$$
\!\!\ba{ll}
\ds
L^{r_1}_\dbF\!(\Om;\!L^{r_2}\!(t,T;\dbR^k)\!)\!=\!\Big\{\f:(t,T)\!\times\!\Om\!\to\!
\dbR^k\; \Big|\; \f(\cd)\hb{
is $\dbF$-adapted and }\dbE\(\!\int_t^T\!\!\!|\f(\tau)|_{\dbR^k}^{r_2}d\tau\!\)^{\frac{r_1}{r_2}}\!\!<\!\infty\!\Big\},\\
\ns\ds
 L^{r_2}_\dbF\!(t,T;\!L^{r_1}\!(\Om;\dbR^k)\!)\!=\!\Big\{\!\f:(t,T)\!\times\!\Om\!\to\!
\dbR^k\; \Big|\; \f(\cd)\hb{ is $\dbF$-adapted and
}\int_t^T\!\!\!\(\!\dbE|\f(\tau)|_{\dbR^k}^{r_1}\!\)^{\frac{r_2}
{r_1}}d\tau\!<\!\infty\!\Big\}.
 \ea
 $$
Both
$L^{r_1}_\dbF(\Om;L^{r_2}(t,T;\dbR^k))$ and
$L^{r_2}_\dbF(t,T;L^{r_1}(\Om;\dbR^k))$ are
Banach spaces with the canonical norms.
In a similar way, we may define
$L^{\infty}_\dbF(\Om;L^{r_2}(t,T;\dbR^k))$,
$L^{r_1}_\dbF(\Om;L^{\infty}(t,T;\dbR^k))$ and
$L^{\infty}_\dbF(\Om;L^{\infty}(t,T;\dbR^k))$.
For $q\in [1,\infty]$, we simply denote
$L^{q}_\dbF(\Om;L^{q}(t,T;\dbR^k))$ by
$L^{q}_\dbF(t,T;\dbR^k)$.  Denote by
$\cS(\dbR^k)$ the set of all $k$-dimensional
symmetric matrices and $I_k$ the $k$-dimensional
identity matrix.

For any $n,m\in\dbN$, and $(s,\eta)\in [0,T)\times
L^2_{\cF_s}(\Om;\dbR^n)$, let us consider the following
controlled linear stochastic differential
equation:
\begin{equation}\label{5.2-eq1}
\left\{ \begin{array}{ll}\ds dx(r)=\big(A(r)
x(r) + B(r) u(r)\big)dr+ \big(C(r) x(r)+D(r)
u(r)\big)dW(r) &\mbox{ in
}[s,T],\\
\ns\ds x(s)=\eta,
\end{array}
\right.
\end{equation}
with the following quadratic cost functional
\begin{equation}\label{5.2-eq2}
\begin{array}{ll}\ds \cJ(s,\eta;u(\cd))
=\frac{1}{2}\mE\Big[ \int_s^T \big(\big\langle Q(r)
x(r),x(r)\big\rangle_{\dbR^n} +\big\langle R(r)
u(r),u(r)\big\rangle_{\dbR^m}\big)dr + \langle
Gx(T),x(T)\rangle_{\dbR^n}\Big].
\end{array}
\end{equation}
In \eqref{5.2-eq1}--\eqref{5.2-eq2}, $
u(\cd)(\in L^2_\dbF(s,T;\dbR^m)$, the space of
admissible controls) is the control variable,
$x(\cd)$ is the state variable, the stochastic
processes $A(\cd)$, $B(\cd)$, $C(\cd)$,
$D(\cd)$, $Q(\cd)$, $R(\cd)$, and the random
variable $G$ satisfy  suitable assumptions to be
given later (see \eqref{20160225e1} in the next
section) such that equation \eqref{5.2-eq1}
admits a unique solution
$x(\cd;s,\eta,u(\cd))\in
L^2_\dbF(\Om;C([s,T];\dbR^n))$, and
\eqref{5.2-eq2} is well defined.  In what
follows, to simplify notations, the time
variable $t$ is sometimes suppressed in $A$,
$B$, $C$, $D$, etc.

In this paper, we are concerned with the
following stochastic linear quadratic control problem (SLQ for short):

\ms

\no\bf Problem (SLQ). \rm For each
$(s,\eta)\in[0,T]\times
L^2_{\cF_s}(\Om;\dbR^n)$, find a $\bar u(\cd)\in
L^2_\dbF(s,T;\dbR^m)$ so that
\begin{equation}\label{5.2-eq3}
\cJ\big(s,\eta;\bar u(\cd)\big)=\inf_{u(\cd)\in
L^2_\dbF(s,T;\dbR^m)}\cJ\big(s,\eta;u(\cd)\big).
\end{equation}

SLQs have been extensively studied in the
literature, for which we refer the readers to \cite{AMZ, Athens, Bismut1, Bismut2,
CLZ1, Tang1, Wonham2, YZ} and the rich references
therein. Similar to the deterministic setting (\cite{Kalman, Wonham1, YL}), Riccati equations (and their variants) are fundamental tools to study SLQs. Nevertheless, for stochastic problems one usually has to consider backward stochastic Riccati equations. For our Problem (SLQ), the desired  backward stochastic Riccati equation takes the following form:
\begin{equation}\label{5.5-eq6}
\left\{
\begin{array}{ll}\ds
dP =-\big( PA + A^{\top}  P + \L C + C^{\top} \L
+ C^{\top}  PC + Q - L^{\top}  K^{\dag} L
\big)dt + \L dW(t) &\mbox{ in }[0,T],\\
\ns\ds P(T)=G,
\end{array}
\right.
\end{equation}
where $A^{\top}$ stands for the transpose of $A$, and
\begin{equation}\label{9.7-eq10}
K\= R+D^{\top}  PD,\qq
L\= B^{\top}  P+D^{\top}  (PC+\L),
\end{equation}
and $K^\dag$ denotes the Moore-Penrose
pseudo-inverse of $K$.

To the authors' best knowledge, \cite{Wonham2} is the first work which employed
Riccati equations to study SLQs.  After
\cite{Wonham2},  Riccati equations were
systematically applied
to study SLQs (e.g. \cite{Athens, Bensoussan1, Bismut2,
Davis, YZ}), and the well-posedness of
such equations was studied in some literatures (See
\cite{YZ, Tang1} and the references cited
therein).

In the early works on SLQs (e.g., \cite{CLZ1, Wonham2, YZ}),  the coefficients $A$, $B$, $C$,
$D$, $Q$, $R$, $G$ appeared in the
control system \eqref{5.2-eq1} and the cost
functional \eqref{5.2-eq2} were assumed to be deterministic. For this case, the corresponding
Riccati equation \eqref{5.5-eq6} is deterministic (i.e., $\Lambda\equiv0$ in \eqref{5.5-eq6}), as well.

To the best of our knowledge, \cite{Bismut1} is the first work that addressed the study of SLQs with random
coefficients. In \cite{Bismut1, Bismut2}, the
author formally derived the equation
\eqref{5.5-eq6}. However, at that time only some
special and simple cases could be solved. Later,
\cite{Peng} proved the well-posedness
for \eqref{5.5-eq6} under the condition that $D=0$ by means of
Bellman's principle of quasi linearization and a
monotone convergence result for symmetric
matrices. This condition was dropped in
\cite{Tang1}, in which it was proved that
\eqref{5.5-eq6} admits a unique solution $(P,\L)$  in a suitable space under the assumptions that  $Q\geq 0$, $G\geq 0$ and $R>\!\!\!>0$.

In Control Theory, one of the fundamental issues
is to find feedback controls, which are
particularly important in practical
applications. It is well-known that, in the
deterministic  case, the purpose to introduce
Riccati equations into the study of Control
Theory (e.g., \cite{Kalman, Wonham1, YL}) is
exactly to design feedback controls for linear
quadratic control problems (LQs for short). More
precisely, under some mild assumptions, one can
show that the unique solvability of
deterministic LQs is equivalent to that of the
corresponding Riccati equations, via which one
can construct the desired optimal feedback
controls. Unfortunately, the same problem is
only partially well-understood in the stochastic
setting, such as the case that all of the
coefficients in \eqref{5.2-eq1}--\eqref{5.2-eq2}
are deterministic (\cite{AMZ, SY}), or the case
that the diffusion term in \eqref{5.2-eq1} is
control-independent, i.e., $D\equiv0$
(\cite{Peng}). However, for the general case, we
shall explain in Remark \ref{20160225r2} below that, the
solution $(P,\L)$ (to \eqref{5.5-eq6}) found in
\cite{Tang1} is not regular enough to serve as
the design of feedback controls for Problem
(SLQ).

Because of the difficulty mentioned above,  it
is quite natural to ask such a question: Is it
possible to link the existence of optimal
feedback controls (rather than the  solvability)
for Problem (SLQ) directly to the solvability of
the equation \eqref{5.5-eq6}? Clearly, from the
viewpoint of applications, it is more desirable
to study the existence of feedback controls for
SLQs than the solvability for the same problems.

The main purpose of this work is to give an affirmative answer to the above question under sharp assumptions on the coefficients appearing in \eqref{5.2-eq1}--\eqref{5.2-eq2}.
For this purpose, let us give the notion of an optimal
feedback operator for
Problem (SLQ).
\begin{definition}\label{5.7-def1}
A stochastic process $\Th(\cd)\in
L^\infty_\dbF(\Om;L^2(0,T;\dbR^{m\times n}))$ is
called an {\it optimal feedback  operator} for
Problem (SLQ) on $[0,T]$ if, for all $(s,\eta)\in
[0,T)\times L^2_{\cF_s}(\Om;\dbR^n)$ and $u(\cd)\in
L^2_\dbF(s,T;\dbR^m)$, it holds that
\begin{equation}\label{5.7-eq2}
\cJ(s,\eta;\Th(\cd)\bar x(\cd))\leq
\cJ(s,\eta;u(\cd)),
\end{equation}
where $\bar x(\cd)=x(\cd\,;s,\eta, \Th(\cd)\bar
x(\cd))$.

\end{definition}

\begin{remark}\label{20160225r1}
In Definition \ref{5.7-def1},
$\Th(\cd)$ is required to be independent of the
initial state $\eta\in L^2_{\cF_s}(\Om;\dbR^n)$.
For a fixed  pair $(s,\eta)\in[0,T)\times
L^2_{\cF_s}(\Om;\dbR^n)$, the inequality
\eqref{5.7-eq2} implies that the control
$$\bar u(\cd)\equiv \Th(\cd)\bar x(\cd)\in L^2_\dbF(s,T;\dbR^m)$$
is optimal for Problem (SLQ). Therefore, for
Problem (SLQ), the existence of an optimal
feedback  operator on $[0,T]$ implies the existence of optimal control for any pair
$(s,\eta)\in [0,T)\times
L^2_{\cF_s}(\Om;\dbR^n)$.
\end{remark}

\begin{remark}\label{20160225r2}
Under some
assumptions, in \cite{Tang1}, it was shown that the equation
\eqref{5.5-eq6} admits a unique solution
$(P,\L)\in L^\infty_\dbF(0,T;\cS(\dbR^n))\times
L^p_\dbF(\Om;L^2(0,T;\cS(\dbR^n)))$  for any
given $p\in [1,\infty)$. Nevertheless  the approach in
\cite{Tang1} does not produce the sharp
regularity $\Th\in L^\infty_\dbF(\Om;L^2(0,T;$
$\cS(\dbR^n)))$ (but rather $\Th\in
L^p_\dbF(\Om;L^2(0,T;\cS(\dbR^n)))$ for any
$p\in [1,\infty)$). Although the author showed
in \cite{Tang1} that if $\bar x$ is an optimal
state, then $\Th\bar x \in L^2_\dbF(0,T;\dbR^n)$ and hence it is
the desired optimal control, such kind of control
strategy is not robust, even with respect to
some very small perturbation. Actually, assume
that there is an error $\d x \in
L^2_\dbF(\Om;C([0,T];\dbR^n))$ (the solution
space of \eqref{5.2-eq1} with $s=0$) with $|\d
x|_{L^2_\dbF(\Om;C([0,T];\dbR^n))}=\e>0$ for
$\e$ being small enough in the observation of
the state, then by the well-posedness result in
\cite{Tang1}, one cannot conclude that $\Th(\bar
x + \d x)$ is an admissible control.
%Even one can do, the norm of $\Th\d x$ may vary large.
Thus, the $\Th$ given in \cite{Tang1} is not a ``qualified" feedback because it is not robust with respect
to small perturbations. How about to assume that $\Th$ has a
good sign or to be monotone (in a suitable sense)? Even for such a special case, it is not hard to see that, things will not become better since we have no other information
about $\d x$ except that it belongs to $L^2_\dbF(\Om;C([0,T];\dbR^n))$, the integrability of the function $\Th\d x$ (with respect to the sample point $\omega$) cannot be improved, and therefore one could not conclude that $\Th(\bar
x + \d x)$ is an admissible control, either.

In a recent paper \cite{Tang2}, the well-posedness result in \cite{Tang1} was slightly improved and it was shown that the solution $(P,\L)$ to
\eqref{5.5-eq6} enjoys the BMO-martingle property. However, this does not help to produce the boundedness of $\Th$ with respect to the sample point $\omega$, either. Actually, we shall give a counterexample (i.e., Example \ref{counterexample-1}) showing that such a boundedness result is not guaranteed without further assumptions.

Let us recall that, the main motivation to
introduce feedback controls is to keep the corresponding
control strategy to be robust with respect to  (small)
perturbations. Hence, the well-posedness results
in \cite{Tang1, Tang2} are not enough to solve our
Problem (SLQ). Nevertheless, for the case that
$D\equiv0$,  the optimal feedback operator in
\eqref{5.10} is specialized as
$$\Th(\cd)=-K(\cd)^{\dag}B(\cd)^\top P(\cd) +
\big(I_m - K(\cd)^{\dag}K(\cd)\big)\th,$$ which
is independent of $\Lambda$, and therefore the
result in \cite{Tang1} (or that in \cite{Peng})
is enough for this special case.
\end{remark}

We have explained that a suitable optimal
feedback control operator for our
Problem (SLQ) should belong to
$L^\infty_\dbF(\Om;L^2(0,T;\dbR^{n\times m}))$.
Nevertheless, to our best knowledge,  the
existence of such operator is completely unknown
for Problem (SLQ) with random coefficients. In
this paper, we shall show
that the existence of the optimal feedback
operator for Problem (SLQ) is equivalent to the
solvability of \eqref{5.5-eq6} in a suitable
sense. When the coefficients $A$, $B$, $C$, $D$,
$G$, $R$, $Q$ are deterministic, such an
equivalence was studied in \cite{AMZ} (see also
\cite{SY} for the problem of a linear quadratic
stochastic two-person zero-sum differential
game). As far as we know, there is no study of
such problems for the general case that $A$, $B$, $C$,
$D$, $R$, $Q$ are stochastic processes, and $G$
is a random variable.

The rest of this paper is organized as follows:
Section \ref{main} is devoted to presenting the
main results of this paper. In Section
\ref{pre}, we give  some preliminary results
which will be used in the remainder of this
paper. Sections \ref{proof}--\ref{proof-2} are
addresses to proofs of our main results. At
last, in Section \ref{example}, we give some
examples for the existence and nonexistence of
the optimal feedback control operator.

%%%%%%%%%%%%%%%%%%%%%%%%%%%%%%%%%%%%%%%%%%%%%%%%

\section{Statement of the main results}\label{main}

%%%%%%%%%%%%%%%%%%%%%%%%%%%%%%%%%%%%%%%%%%%%%%%%

Let us first introduce the following assumption:

\ms

\no({\bf AS1}) {\it The coefficients in \eqref{5.2-eq1}--\eqref{5.2-eq2} satisfy the following measurability/integrability conditions:}
\bel{20160225e1}
\begin{array}{ll}\ds
A(\cd)\in
L^\infty_\dbF(\Om;L^1(0,T;\dbR^{n\times n})), \q
 C(\cd)\in
L^\infty_\dbF(\Om;L^2(0,T;\dbR^{n\times n})), \\
B(\cd)\in
L^\infty_\dbF(\Om;L^2(0,T;\dbR^{n\times m})), \q
D(\cd)\in L^\infty_\dbF(0,T;\dbR^{n\times m}),\\
\ns\ds Q(\cd)\in L^\infty_\dbF(\Omega;L^1(0,T;\cS(\dbR^n))),
\q R(\cd) \in L^\infty_\dbF(0,T;\cS(\dbR^m)),\q
G\in L^\infty_{\cF_T}(\Om;\cS(\dbR^n)).
\end{array}
\ee

We have the following result:

\begin{theorem}\label{5.7-th1}
Let the assumption (AS1) hold. Then, Problem (SLQ) admits an
optimal feedback operator $\Th(\cd)\in
L^\infty_\dbF(\Om;L^2(0,T;\dbR^{m\times n}))$ if
and only if the Riccati equation \eqref{5.5-eq6}
admits a solution
$\big(P(\cd),\L(\cd)\big)\in
L^{\infty}_{\dbF}(\Om;C([0,T];\cS(\dbR^n)))
\times L^p_{\dbF}(\Omega;L^2(0,T;\cS(\dbR^n)))$
(for all $p\geq 1$)
such that
\bel{20160225e2}
 \cR(K(t,\omega))\supset\cR(L(t,\omega)) \q\hbox{and}\q K(t,\omega)\geq 0,\qq \ae (t,\omega)\in
 [0,T]\times\Omega,
 \ee
and
\begin{equation}\label{5.7-eq5}
K(\cd)^{\dag}L(\cd)\in
L^\infty_\dbF(\Om;L^2(0,T;\dbR^{m\times n})).
\end{equation}
In this case, the optimal feedback operator
$\Th(\cd)$ is given as
\begin{equation}\label{5.10}
\Th(\cd)=-K(\cd)^{\dag}L(\cd) + \big(I_m -
K(\cd)^{\dag}K(\cd)\big)\th,
\end{equation}
where $\th\in
L^\infty_\dbF(\Om;L^2(0,T;\dbR^{m\times n}))$ is
arbitrarily given. Furthermore,
\begin{equation}\label{Value}
\inf_{u\in
L^2_\dbF(s,T;\dbR^m)}\cJ(s,\eta;u)=\frac{1}{2}\,\dbE\langle
P(s)\eta,\eta\rangle_{\dbR^n}.
\end{equation}
\end{theorem}

\begin{corollary}\label{corollary-1}
Let (AS1) hold. Then the Riccati equation
\eqref{5.5-eq6} admits at most one solution
$\big(P(\cd),\L(\cd)\big)\in
L^{\infty}_{\dbF}(\Om;C([0,T];\cS(\dbR^n)))
\times L^p_{\dbF}(\Om;L^2(0,T;\cS(\dbR^n)))$
(for all $p\geq 1$) satisfying
\eqref{20160225e2} and \eqref{5.7-eq5}.
\end{corollary}

%Noting that we do not assume the nonnegative/positive
%definiteness of $R$, $Q$ and $G$, the solution
%to Problem (SLQ) may not exist. In this paper,
%since we focus on the existence of the optimal
%feedback control rather than the existence of
%the optimal control, let us make further the following
%assumption.

%\ms

%\no({\bf AS2}) {\it For any $(s,\eta)\in[0,T)\times
%L^2_{\cF_s}(\Om;H)$, there exists a unique
%optimal control for Problem (SLQ).}

%\ms

%It is easy to see that, when $Q\geq 0$, $G\geq 0$ and $R>\!\!\!>0$, the assumption (AS2) is satisfied.

%Under the additional assumption (AS2),

The result in Theorem \ref{5.7-th1} can be strengthened as follows.

\begin{theorem}\label{5.7-th1-2}
Let the assumption (AS1) hold. Then, Problem (SLQ) admits a unique
optimal feedback operator $\Th(\cd)\in
L^\infty_\dbF(\Om;L^2(0,T;\dbR^{m\times n}))$ if
and only if the Riccati equation \eqref{5.5-eq6}
admits a unique solution
$\big(P(\cd),\L(\cd)\big)\in
L^{\infty}_{\dbF}(\Om;C([0,T];\cS(\dbR^n)))
\times L^p_{\dbF}(\Om;L^2(0,T;$ $\cS(\dbR^n)))$
(for all $p\geq 1$)
such
that  $K(t,\omega)> 0$ for $\ae (t,\omega)\in
 [0,T]\times\Omega$ (and hence $ K^{\dag}$ in \eqref{5.5-eq6} can be replaced by $ K^{-1}$) and
$
K(\cd)^{-1}L(\cd)\in
L^\infty_\dbF(\Om;L^2(0,T;\dbR^{m\times n}))$.
In this case, the optimal feedback operator
$\Th(\cd)$ is given by $\Th(\cd)=-K(\cd)^{-1}L(\cd)$, and \eqref{Value} (in Theorem \ref{5.7-th1}) holds.
\end{theorem}

Several remarks are in order.

\begin{remark}
We borrow some idea from \cite{AMZ, SY} to employ the Moore-Penrose
pseudo-inverse in the study of Riccati
equations for SLQs when the matrix $K$ in \eqref{9.7-eq10} is singular.
\end{remark}

\begin{remark}
The proof of sufficiency in Theorems
\ref{5.7-th1}--\ref{5.7-th1-2} is very close to
the deterministic setting and also that of the
case that the coefficients in
\eqref{5.2-eq1}--\eqref{5.2-eq2} are
deterministic. The main difficulty in the proof
of necessity in Theorems
\ref{5.7-th1}--\ref{5.7-th1-2} consists in the
very fact that the equation \eqref{5.5-eq6} is a
nonlinear equation with a non-global Lipschitz
nonlinearity. Nevertheless, since Riccati
equations appearing in Control Theory enjoy some
special structures, at least under some
assumptions they are still globally solvable. A
basic idea to solve Riccati equations globally is to link
them with suitable solvable optimal control
problems, and via which one obtains the desired
solutions. To the best of our knowledge, such an
idea was first used to solve deterministic
differential Riccati equations in \cite{Reid}
(though in that paper, the author considered the
second variation for a nonsingular nonparametric
fixed endpoint problem in the calculus of
variations rather than an optimal control
problem). This idea was later adopted by  many
authors (e.g., \cite{AMZ, Bismut2, Kalman, SY,
Tang1}). In this work, we shall also use such an
idea.
\end{remark}

\begin{remark}
To simplify the presentation, in this paper we assume that the filtration $\dbF$ is natural. One can also consider the case of general filtration. Of course, for general filtration the solutions to \eqref{5.5-eq6} have to be understood in the sense of transposition (introduced in \cite{LZ0, LZ1}).
\end{remark}

\begin{remark}
The same SLQ problems (as those in this paper) but in infinite dimensions still make sense. However, the new difficulty in the infinite dimensional setting is how to explain the stochastic integral $\int_0^T\L (t)dW(t)$ that appeared in \eqref{5.5-eq6} because for this case $\L(\cd)$ is an operator-valued stochastic process, and therefore one has to use the theory of transposition solution for operator-valued backward stochastic evolution equations (\cite{LZ1, LZ}). Progress in this respect is presented in \cite{LZ3}.
\end{remark}

\begin{remark}
It would be quite interesting to extend the main result in this paper to linear quadratic
stochastic differential games or similar problems for mean-field stochastic differential equations. Some relevant studies can be found in \cite{Huyen, SY} but the full pictures are still unclear.
\end{remark}

%%%%%%%%%%%%%%%%%%%%%%%%%%%%%%%%%%%%%%%%%%%%%%%%

\section{Some preliminary results}\label{pre}

%%%%%%%%%%%%%%%%%%%%%%%%%%%%%%%%%%%%%%%%%%%%%%%%

In this section, we present some preliminary
results, which will be useful later.

First, for any $s\in[0,T)$, we consider the
following stochastic differential equation:
\begin{equation}\label{6.20-eq1}
\left\{
\begin{array}{ll}\ds
dx = (\cA x  + f)dt + (\cB x+g)dW(t) &\mbox{ in }[s,T],\\
\ns\ds x(s)=\eta.
\end{array}
\right.
\end{equation}
Here $\cA,\cB\in
L^\infty_\dbF(\Om;L^2(0,T;\dbR^{k\times k}))$,
$\eta\in L^2_{\cF_s}(\Om;\dbR^k)$, and $f,g\in
L^2_\dbF(s,T;\dbR^k)$.

Let us recall the following result (We refer to \cite[Chapter V, Section
3]{Protter} for its proof).

\begin{lemma}\label{lm2}
The equation \eqref{6.20-eq1} admits one and only one $\dbF$-adapted solution $x(\cd)\in
L^2_\dbF(\Om;$ $C([s,T];\dbR^k))$.
\end{lemma}

Next, we need to consider the following backward stochastic
differential equation:
\begin{equation}\label{2.14-eq1}
\left\{
\begin{array}{ll}\ds
dy = f(t,y,z)dt + zdW(t) &\mbox{ in }[s,T],\\
\ns\ds y(T)=\xi.
\end{array}
\right.
\end{equation}
Here  $\xi\in L^\infty_{\cF_T}(\Om;\dbR^k)$, and
$f$ satisfies that
\bel{20160226e5}
\begin{cases}
f(\cd,0,0)\in
L^\infty_\dbF(\Om;L^1(s,T;\dbR^k)),\\
|f(\cd,\a_1,\a_2)-f(\cd,\b_1,\b_2)|_{\dbR^k}\leq
f_1(\cd)|\a_1-\b_1|_{\dbR^k} +
f_2(\cd)|\a_2-\b_2|_{\dbR^k}, \;\forall
\a_1,\a_2,\b_1,\b_2\in \dbR^k,
\end{cases}
\ee
where $f_1(\cd)\in L^\infty_\dbF(\Om;L^1(s,T;\dbR))$
and $f_2(\cd)\in L^\infty_\dbF(\Om;L^2(s,T;\dbR))$.

By means of \cite[Theorem 2.7]{Delbaen-Tang} (See also \cite{BDHPS} for an early result in this direction), we have
\begin{lemma}\label{lm2.1}
For any $p>1$, the equation \eqref{2.14-eq1} admits one and only one $\dbF$-adapted solution
$(y(\cd),z(\cd))\in
L^\infty_\dbF(\Om;C([s,T];\dbR^k))\times
L^p_\dbF(\Om;L^2(s,T;\dbR^k))$.
%
%with any $p>1$, and
%%
%\begin{equation}\label{lm2.1-eq1}
%|(y(\cd),z(\cd))|_{L^\infty_\dbF(\Om;C([0,T];\dbR^d))\times
%L^p_\dbF(\Om;L^2(0,T;\dbR^d))}\leq
%\cC\big(|\xi|_{L^\infty_{\cF_T}(\Om;\dbR^d)} +
%|f(\cd,0,0)|_{L^\infty_\dbF(\Om;L^1(0,T;\dbR^d))}
%\big).
%\end{equation}
%
\end{lemma}

\ms

Further, let us recall the following known
Pontryagin-type maximum principle (\cite[Theorem 3.2]{Bismut1}).

\begin{lemma}\label{5.4-th1}
Let $(\bar x(\cd),\bar u(\cd))$ $\in
L^2_\dbF(\Om;C([s,T];\dbR^n))\times
L^2_\dbF(s,T;\dbR^m)$ be an optimal pair of
Problem (SLQ). Then there exists a pair
$ (\bar y(\cd), \bar z(\cd)) \in
L^2_\dbF(\Om;C([s,T];\dbR^n))\times
L^2_\dbF(s,T;\dbR^n) $
satisfying the following backward stochastic
evolution equation:
$$
\left\{
\begin{array}{ll}\ds
d\bar y(t)=-\big(A^{\top} \bar y(t)+C^{\top}
\bar z(t)+Q\bar x(t)\big)dt +\bar z(t)dW(t)
&\mbox{\rm in }[s,T],\cr
\ns\ds\bar y(T)= G\bar x(T),
\end{array}
\right.
$$
and
$$
R\bar u(\cd)+B^{\top} \bar y(\cd)+D^{\top} \bar z(\cd) =0,\qq\ae (t,\om)\in [s,T]\times\Om.
$$
\end{lemma}

As an immediate consequence of Lemmas \ref{lm2}
and \ref{5.4-th1}, we have the following result.
\begin{corollary}\label{5.7-prop1}
Let $\Th(\cd)$ be an optimal feedback operator
for Problem (SLQ). Then, for any
$(s,\eta)\in[0,T)\times
L^2_{\cF_s}(\Om;\dbR^n)$, the following
forward-backward stochastic differential
equation:
$$
\left\{
\begin{array}{ll}
\ns\ds d\bar x(t) =(A+B\Th)\bar x(t) dt+
(C+D\Th)\bar x(t) dW(t)
&\hb{\rm in } [s,T],\\
\ns\ds
d\bar y(t)=-\big(A^{\top}  \bar y(t)+C^{\top}  \bar z(t)+ Q\bar x(t)\big)dt+\bar z(t)dW(t)\quad&\mbox{\rm in } [s,T],\\
\ns\ds \bar x(s)=\eta,\q\bar y(T)=G\bar x(T),
\end{array}
\right.
$$
admits a unique solution
$ (\bar x(\cd),\bar y(\cd),\bar z(\cd))\in
L^2_\dbF(\Om;C([s,T];\dbR^n))\times
L^2_\dbF(\Om;C([s,T];\dbR^n))\times
L^2_\dbF(s,T;\dbR^n) $,
and
$$
R\Th \bar x(\cd)+B^{\top}  \bar y(\cd) +D^{\top}  \bar z(\cd)
=0,\q \ae (t,\om)\in [s,T]\times\Om.
$$
\end{corollary}

\ms

Finally, for the reader's convenience, let us
recall the following result for the
Moore-Penrose pseudo-inverse and refer the
readers to \cite[Chapter 1]{Ben-Israel} for its
proof.

\begin{lemma}\label{lm3}
{\rm 1)} Let $M\in\dbR^{n\times n}$. Then the
Moore-Penrose pseudo-inverse $M^\dag$ of $M$
satisfies that
$$
M^\dag = \lim_{\delta \searrow 0} (M^\top M +
\delta I_n)^{-1} M^\top.
$$

{\rm 2)}  If $M\in \cS(\dbR^n)$, then $M^\dag M = M
M^\dag$ and $M^\dag M$ is the orthogonal
projector from $\dbR^n$ to the range of $M$.
\end{lemma}
%

%%%%%%%%%%%%%%%%%%%%%%%%%%%%%%%%%%%%%%%%%%%%%%%%

\section{Proof of  Theorem
\ref{5.7-th1}}\label{proof}

%%%%%%%%%%%%%%%%%%%%%%%%%%%%%%%%%%%%%%%%%%%%%%%%

 In this section, we give a proof of Theorem
\ref{5.7-th1}.

\subsection{Proof of sufficiency in Theorem \ref{5.7-th1}}

In this subsection, we prove the ``if" part in Theorem \ref{5.7-th1}. The proof is more or less standard. For the reader's convenience, we provide here the details.

Let us assume that equation \eqref{5.5-eq6} admits
a solution
$$\big(P(\cd),\L(\cd)\big)\in
L^{\infty}_{\dbF}(\Om;C([0,T];\cS(\dbR^{n})))
\times
L^p(\Om;L^2_{\dbF}(0,T;\cS(\dbR^{n})))$$
 so
that \eqref{20160225e2} and \eqref{5.7-eq5} hold. Then, for any $\th\in
L^\infty_\dbF(\Om;$ $L^2(0,T;\dbR^{m\times n}))$, by \eqref{9.7-eq10} and  \eqref{5.7-eq5}, the function
$ \Th(\cd)$ given by \eqref{5.10} belongs to $L^\infty_{\dbF}(\Om;L^2(0,T;$ $\dbR^{m\times n}))$.
For any $s\in [0,T)$, $\eta\in
L^2_{\cF_s}(\Om;\dbR^n)$, and $u(\cd)\in
L^2_\dbF(s,T;\dbR^m)$, let $x(\cd)\equiv
x(\cd\,;s,\eta,u(\cd))$ be the corresponding
state process for \eqref{5.2-eq1}.
By It\^o's formula, and using \eqref{5.2-eq1}, \eqref{5.5-eq6}, \eqref{9.7-eq10}, we obtain that
\begin{equation}\label{Ito-formula-to-P}
\begin{array}{ll}\ds
 d \big\langle P
x ,x \big\rangle_{\dbR^n}\\
\ns\ds= \big\langle d P
x ,x \big\rangle_{\dbR^n} + \big\langle
P dx ,x \big\rangle_{\dbR^n} +
\big\langle P
x ,dx \big\rangle_{\dbR^n}\\
\ns\ds\q + \big\langle d
P
dx ,x \big\rangle_{\dbR^n}+ \big\langle dP
x ,dx \big\rangle_{\dbR^n} + \big\langle
P  dx ,dx \big\rangle_{\dbR^n} \\
\ns\ds =\big\langle -
\big[PA+ A^\top P+\L C
+C^\top \L+C^\top P C+ Q - L^{\top}  K^{\dag} L
\big]x,x \big\rangle_{\dbR^n}dr\\
\ns\ds \q  + \big\langle P (A x  +
Bu) ,x\big\rangle_{\dbR^n}dr+
\big\langle P (C
x  + Du) ,x\big\rangle_{\dbR^n}dW(r)\\
\ns\ds \q + \big\langle Px,A x +
Bu\big\rangle_{\dbR^n}dr + \big\langle
Px
,C x  + Du\big\rangle_{\dbR^n}dW(r)\\
\ns\ds \q + \big\langle \L (C x  +
Du),x\big\rangle_{\dbR^n}dr +
\big\langle \L x,C x  +
Du\big\rangle_{\dbR^n}dr\\
\ns\ds\q + \big\langle P (C x  +
Du),C x  +
Du\big\rangle_{\dbR^n}dr + \big\langle
\L x, x\big\rangle_{\dbR^n}dW(r)\\
\ns\ds =  -\big\langle( Q-L^\top K^{\dag}L)x,x \big\rangle_{\dbR^n}dr
+\big\langle P Bu,
x\big\rangle_{\dbR^n}dr
 \\
\ns\ds \q + \big\langle P x,
Bu\big\rangle_{\dbR^n}dr+ \big\langle P C
x, Du\big\rangle_{\dbR^n}dr +
\big\langle PDu, Cx+
Du\big\rangle_{\dbR^n}dr\\
\ns\ds \q + \big\langle Du,
\L x\big\rangle_{\dbR^n}dr + \big\langle
\L x, Du\big\rangle_{\dbR^n}dr +
\langle P (C x  + Du)
,x\rangle_{\dbR^n}dW(r)\\
\ns\ds \q + \langle Px ,C x  +
Du\rangle_{\dbR^n}dW(r)+ \big\langle \L
x, x\big\rangle_{\dbR^n}dW(r)\\
\ns\ds =-\big\langle( Q-L^\top K^{\dag}L)x,x \big\rangle_{\dbR^n}dr+2\langle L^{\top}u,x\rangle_{\dbR^n}dr+\langle D^{\top}P Du,u\rangle_{\dbR^m} dr\\
\ns\ds\q +\big[2\langle P(Cx+Du),x\rangle_{\dbR^n}+\langle \Lambda x,x\rangle_{\dbR^n}\big]dW(r).
\end{array}
\end{equation}

Since $K$ is an adapted process, from the first conclusion
in Lemma \ref{lm3}, we deduce that $K^\dag$ is also
adapted.

Notice that from \eqref{5.10} one has
$$\ba{ll}
\ns\ds K \Th=-KK^{\dagger}L,\qq L+K \Th
=L-KK^{\dagger}L.
\ea$$
Moreover, by $\cR(K(\cd))\supset\cR(L(\cd))$, we
conclude that for a.e. $(t,\om)\in (0,T)\times\Om$,
and for any $v\in \dbR^n$, there is a $\hat
v\in\dbR^n$ such that $K(t,\om)\hat v = L(t,\om)
v$. Hence
$$\ba{ll}
\ds  L(t,\om)v+K(t,\om) \Th(t,\om)v
\\
\ns\ds=L(t,\om)v-K(t,\om)
K^{\dagger}(t,\om)L(t,\om)v \\
\ns\ds =K
(t,\om)v-K(t,\om)K^{\dagger}(t,\om)K(t,\om)\hat
v=0.
\ea$$
This yields that
$$
L(t,\om)v+K(t,\om) \Th(t,\om)=0 \;\mbox{ for
a.e. }  (t,\om)\in (0,T)\times\Om,
$$
which, together with the symmetry of $K(\cd)$,
implies that $L^{\top}=-\Th^{\top}K$. Since in
this case $\Th(\cd)\in
L^\infty_{\dbF}(\Omega;L^2(0,T;\dbR^{m\times
n}))$, $K(\cd)$ is bounded, one has $L(\cd)\in
L^\infty_{\dbF}(\Omega;L^2(0,T;\dbR^{m\times
n}))$. Moreover, from the definition of $\Th$ in
(\ref{5.10}), we derive that,
$$
\Th^{\top}K\Th=-\Th^{\top}K\big[K^{\dagger}L+(I_m-K^{\dagger}K)\th\big]=-\Th^{\top}K K^{\dagger}L=L^{\top}K^{\dagger} L.
$$
As a result, we rewrite (\ref{Ito-formula-to-P}) as,
\bel{Ito-formula-to-P-2}
 \ba{ll}
 \ns\ds d \big\langle P
x ,x \big\rangle_{\dbR^n} =-\big\langle( Q-\Th^\top K\Th)x,x \big\rangle_{\dbR^n}dr+2\langle L^{\top}u,x\rangle_{\dbR^n}dr+\langle D^{\top}P Du,u\rangle_{\dbR^m} dr\\
\ns\ds\qq\qq\qq+\big[2\langle P(Cx+Du),x\rangle_{\dbR^n}+\langle \Lambda x,x\rangle_{\dbR^n}\big]dW(r).
\ea
\ee
In order to deal with the stochastic integral
above, for any $s\in[0,T),$ we introduce the
following sequence of stopping times $\t_j$ as,
$$
\t_j \triangleq \inf\Big\{t\geq s\;\Big|\;
\int_s^t|\L(r)|^2dr\geq j\Big\}\wedge T.
$$
It is easy to see that $\t_j\rightarrow T$,
$\dbP$-a.s., as $j\rightarrow\infty.$ Using
(\ref{Ito-formula-to-P-2}), we obtain that,
$$
\begin{array}{ll}
\ns\ds \dbE\langle
P(\t_j)x(\t_j),x(\t_j)\rangle_{\dbR^n}
+\dbE\int_s^T\chi_{[s,\tau_j]}\big[\langle Q
x(r),x(r)\rangle_{\dbR^n}
+\langle Ru(r),u(r)\rangle_{\dbR^m}\big]dr\\
\ns\ds=\dbE\langle
P(s)\eta,\eta\rangle_{\dbR^n}+
\dbE\int_s^T\chi_{[s,\tau_j]}\big[\langle\Th^{\top}K\Th
x(r),x(r)\rangle_{\dbR^n}
+2 \langle L^{\top}u(r),x(r)\rangle_{\dbR^n}\big]dr\\
\ns\ds\q +\dbE\int_s^T \chi_{[s,\tau_j]}\langle
(R+D^{\top}P D)u(r),u(r)\rangle_{\dbR^m} dr.
\end{array}
$$
Clearly,
$$
|\langle
P(\t_j)x(\t_j),x(\t_j)\rangle_{\dbR^n}|\leq
|P|_{L^\infty_\dbF(0,T;\dbR^{n\times
n})}|x|^2_{L^2_\dbF(\Om;C([0,T];\dbR^{n})},
$$
by Dominated Convergence Theorem, we obtain
that
\begin{equation}\label{4.7-eq1}
\lim_{j\to\infty}\langle
P(\t_j)x(\t_j),x(\t_j)\rangle_{\dbR^n}= \langle
P(T)x(T),x(T)\rangle_{\dbR^n}.
\end{equation}
Furthermore,
$$
\begin{array}{ll}\ds
\big|\chi_{[s,\tau_j]}\big[\langle Q
x(r),x(r)\rangle_{\dbR^n} +\langle
Ru(r),u(r)\rangle_{\dbR^m}\big]\big|\\
\ns\ds \leq \big|\big[\langle Q
x(r),x(r)\rangle_{\dbR^n} +\langle
Ru(r),u(r)\rangle_{\dbR^m}\big]\big|\in
L^1_\dbF(0,T),
\end{array}
$$
by Dominated Convergence Theorem again, we
obtain that
\begin{equation}\label{4.7-eq2}
\begin{array}{ll}\ds
\lim_{j\to\infty}\dbE\int_s^T\chi_{[s,\tau_j]}\big[\langle
Q x(r),x(r)\rangle_{\dbR^n} +\langle
Ru(r),u(r)\rangle_{\dbR^m}\big]dr \\
\ns\ds =\dbE\int_s^T \big[\langle Q
x(r),x(r)\rangle_{\dbR^n} +\langle
Ru(r),u(r)\rangle_{\dbR^m}\big]dr.
\end{array}
\end{equation}
Similarly, we show that
\begin{equation}\label{4.7-eq3}
\begin{array}{ll}
 \ds \lim_{j\to\infty}
\dbE\int_s^T\chi_{[s,\tau_j]}\big[\langle\Th^{\top}K\Th
x(r),x(r)\rangle_{\dbR^n}
+2 \langle L^{\top}u(r),x(r)\rangle_{\dbR^n}\big]dr\\
\ns\ds\q +\lim_{j\to\infty} \dbE\int_s^T
\chi_{[s,\tau_j]}\langle (R+D^{\top}P
D)u(r),u(r)\rangle_{\dbR^m} dr \\
\ns\ds = \dbE\int_s^T \big[\langle\Th^{\top}K\Th
x(r),x(r)\rangle_{\dbR^n}
+2 \langle L^{\top}u(r),x(r)\rangle_{\dbR^n}\big]dr\\
\ns\ds\q + \dbE\int_s^T \langle (R+D^{\top}P
D)u(r),u(r)\rangle_{\dbR^m} dr.
\end{array}
\end{equation}

It follows from \eqref{4.7-eq1}--\eqref{4.7-eq3}
that
\begin{equation}\label{5.31-eq1}
\begin{array}{ll}
\ds 2\cJ(s,\eta;u(\cd))\\
 \ns\ds=
\dbE\langle Gx(T),x(T)\rangle_{\dbR^n}+\dbE\int_s^{T}\big[\langle Q x(r),x(r)\rangle_{\dbR^n}
+\langle Ru(r),u(r)\rangle_{\dbR^m}\big]dr\\
\ns\ds=\dbE\langle P(s)\eta,\eta\rangle_{\dbR^n}+
\dbE\int_s^{T}\big[\langle\Th^{\top}K\Th x ,x\rangle_{\dbR^n}
+2 \langle L^{\top}u ,x \rangle_{\dbR^n}+\langle K u ,u\rangle_{\dbR^m}
\big]dr\\
\ns\ds =\dbE\Big[\big\langle
P(s)\eta,\eta\big\rangle_{\dbR^n}+\int_s^T\big(\big\langle
K\Th x,\Th x\big\rangle_{\dbR^m}-2\big\langle K\Th x,u\big\rangle_{\dbR^m}+\langle Ku,u\rangle_{\dbR^m}\big)dr\Big]\\
\ns\ds =2\cJ(s,\eta;\Th\bar x)
+\dbE\int_s^T\big\langle K(u-\Th x),u-\Th
x\big\rangle_{\dbR^m}dr,
\end{array}
\end{equation}
where we have used the fact that $L^{\top}=-\Th^{\top}K$. Hence, by $K(\cd)\geq0$, we have
$$
\cJ(s,\eta;\Th\bar x)\leq \cJ(s,\eta;u),\q\forall\,
u(\cd)\in L^2_\dbF(s,T;\dbR^m).
$$
Thus, for any $\th\in L^\infty_\dbF(\Om;$
$L^2(0,T;\dbR^{m\times n}))$, the function $
\Th(\cd)$ given by \eqref{5.10} is an optimal
feedback operator for Problem (SLQ). This
completes the proof of sufficiency in
Theorem \ref{5.7-th1}.

\subsection{Proof of necessity in Theorem \ref{5.7-th1}}

This subsection is addressed to proving the ``only if" part in Theorem \ref{5.7-th1}.
 We borrow some ideas from
\cite{AMZ, Bismut2, Kalman, Reid, SY}, and divide the proof into several steps.

\ms

{\bf Step 1}. Let $\Th(\cd)\in
L^\infty_\dbF(\Om;L^2(0,T;\dbR^{m\times n}))$ be
an optimal feedback operator for Problem (SLQ)
on $[0,T]$. Then, by Corollary \ref{5.7-prop1},
for any $\zeta\in \dbR^n$, the following
forward-backward stochastic differential
equation
\begin{equation}\label{5.7-eq4.1}
\left\{
\begin{array}{ll}
\ds dx(t)= (A+B\Th)x(t)dt+ (C+D\Th)x(t)
dW(t)& \mbox{ in }[0,T],\\
\ns\ds
dy(t)=-\big(A^{\top}  y(t)+C^{\top}  z(t)+Qx(t)\big)dt+z(t)dW(t) & \mbox{ in }[0,T],\\
\ns\ds x(0)=\zeta,\qq y(T)= Gx(T)
\end{array}
\right.
\end{equation}
admits a solution $(x(\cd),y(\cd),z(\cd))\in
L^2_\dbF(\Om;C([0,T];\dbR^n))\times
L^2_\dbF(\Om;C([0,T];\dbR^n))\times
L^2_\dbF(0,T;\dbR^n)$ so that
\begin{equation}\label{5.7-eq3}
R\Th x+B^{\top}  y+D^{\top}  z =0,\q\ae
(t,\om)\in (0,T)\times\Om.
\end{equation}

Also, consider the following stochastic
differential equation:
\begin{equation}\label{5.26-eq4}
\left\{
\begin{array}{ll}
\ns\ds d\tilde x =
\big[-A-B\Th+\big(C+D\Th\big)^2 \big]^{\top}
\tilde x dt - \big(C+D\Th\big)^{\top} \tilde
x dW(t) \qq  \mbox{ in }[0,T],\\
\ns\ds \tilde x(0)=\zeta.
\end{array}
\right.
\end{equation}
By Lemma \ref{lm2}, the equation
\eqref{5.26-eq4} admits a unique solution
$\tilde x\in L^2_\dbF(\Om;C([0,T];\dbR^n))$.

Further, consider the following $\dbR^{n\times
n}$-valued equations:
\begin{equation}\label{7.19-eq1}
\left\{
\begin{array}{ll}
\ds dX=  (A+B\Th)X dt+ (C+D\Th) X
dW(t)& \mbox{ in }[0,T],\\
\ns\ds dY =-\big(A^{\top} Y+ C^{\top} Z+
QX\big)dt+ ZdW(t) &
\mbox{ in }[0,T],\\
\ns\ds X(0)=I_{n}, \q Y(T)= GX(T)
\end{array}
\right.
\end{equation}
and
\begin{equation}\label{7.19-eq2}
\left\{
\begin{array}{ll}
\ds d\wt X= \big[-A-B\Th+\big(C+D\Th\big)^2
\big]^{\top} \wt X dt - (C+D\Th)^{\top} \wt X
dW(t)& \mbox{ in }[0,T],\\
\ns\ds \wt X(0)=I_n.
\end{array}
\right.
\end{equation}
In view of Corollary \ref{lm2}, it is easy to
show that equations \eqref{7.19-eq1} and
\eqref{7.19-eq2} admit, respectively, unique
solutions $ (X,Y,Z)\in
L^2_\dbF(\Om;C([0,T];\dbR^{n\times n}))\times
L^2_\dbF(\Om;C([0,T];\dbR^{n\times n}))\times
L^2_\dbF(0,T;\dbR^{n\times n}) $ and $ \wt X \in
L^2_\dbF(\Om;C([0,T];\dbR^{n\times n}))$.

It follows from \eqref{5.7-eq4.1} to
\eqref{7.19-eq2}  that, for any $\zeta\in \dbR^n$,
\begin{equation}\label{7.19-eq6}
\begin{array}{ll}\ds
x(t;\zeta) =  X(t)\zeta, \q
y(t;\zeta) =  Y(t)\zeta, \q \tilde
x(t;\zeta) = \wt X(t)\zeta,\q& \forall\, t\in [0,T],\\[1mm]
z(t;\zeta) =  Z(t)\zeta, & \ae t\in [0,T].
\end{array}
\end{equation}

By \eqref{5.7-eq3} and noting \eqref{7.19-eq6},
we find that
\begin{equation}\label{5.7-eq5.1}
R\Th X+B^{\top}  Y+D^{\top}  Z =0, \q \ae
(t,\om)\in[0,T]\times\Om.
\end{equation}

For any $\zeta,\rho\in \dbR^n$ and $t\in [0,T]$, by It\^o's formula, we have
$$
\begin{array}{ll}\ds\big\langle
x(t;\zeta),\tilde
x(t;\rho)\big\rangle_{\dbR^n} -
\big\langle \zeta,\rho\big\rangle_{\dbR^n}\\
\ns\ds =\int_0^t \big\langle  \big(A +B
\Th\big) x(r;\zeta),\tilde
x(r;\rho)
\big\rangle_{\dbR^n} dr + \int_0^t \big\langle \big(C +D
\Th \big)
x(r;\zeta),\tilde x(r;\rho) \big\rangle_{\dbR^n} dW(r)\\
\ns\ds \q + \int_0^t \big\langle
x(r;\zeta), \big[-A -B \Th
+\big(C+D\Th\big)^2\big]^* \tilde
x(r;\rho)
\big\rangle_{\dbR^n} dr\\
\ns\ds \q - \int_0^t \big\langle
x(r;\zeta), \big(C+D\Th\big)^*
\tilde
x(r;\rho) \big\rangle_{\dbR^n} dW(r) \\
\ns\ds \q- \int_0^t \big\langle \big(C +D
\Th \big) x(r;\zeta), \big(C +D
\Th \big)^* \tilde x(r;\rho)
\big\rangle_{\dbR^n} dr\\
\ns\ds =0.
\end{array}
$$
Thus,
$$
\big\langle X(t)\zeta, \wt
X(t)\rho\big\rangle_{\dbR^n}=\big\langle
x(t;\zeta),\tilde
x(t;\rho)\big\rangle_{\dbR^n} = \big\langle
\zeta,\rho\big\rangle_{\dbR^n}, \q \dbP\mbox{-}\as
$$
This implies that $X(t)\wt X(t)^*=I_n$,
$\dbP$-a.s., that is, $\wt
X(t)^*=X(t)^{-1}$, $\dbP$-a.s.

\medskip

{\bf Step 2}. Put
\begin{equation}\label{8.25-eq3}
P(t,\om)\triangleq Y(t,\om) \wt X(t,\om)^{\top}, \q
\Pi(t,\om)\triangleq Z(t,\om) \wt X(t,\om)^{\top}.
\end{equation}
By It\^o's formula,
$$
\3n\begin{array}{ll} \ds
dP\3n&\ds=\Big\{-\big(A^{\top}Y+ C^{\top} Z+
Q X \big)X^{-1} +  Y X^{-1} \big[(C
+D \Th )^2  - A -B \Th \big]\\
\ns&\ds\q- Z X^{-1} (C +D \Th ) \Big\}dt
+\[ Z X^{-1} - Y X^{-1} (C +D \Th )\]dW(t)\\
\ns&\ds =\Big\{- A^{\top} P - C^{\top} \Pi - Q
+P \big[(C+D\Th)^2-A-B\Th\big]-\Pi (C +D\Th) \Big\}dt \\
\ns&\ds \q +\big[\Pi-P(C+D\Th) \big]dW(t).
\end{array}
$$
Let
\begin{equation}\label{9.7-eq11}
\L\triangleq\Pi-P (C+D\Th).
\end{equation}
Then, $(P(\cd),\L(\cd))$ solves the following
$\dbR^{n\times n}$-valued backward stochastic
differential equation:
\begin{equation}\label{6.7-eq8}
\left\{
\begin{array}{ll}\ds
dP =-\big[PA+ A^{\top}P + \L C+ C^{\top}\L+
C^{\top} P C \\
\ns\ds \qq\q +(P B + C^{\top}PD + \L D) \Th + Q
\big]dt + \L dW(t) &\mbox{ in }[0,T],\cr
\ns\ds P(T)=G.
\end{array}
\right.
\end{equation}
By Lemma \ref{lm2.1}, we conclude that
$(P,\L)\in
L^\infty_\dbF(\Om;C([0,T];\dbR^{n\times n}))
\times L^p_\dbF(\Om;L^2(0,T;\dbR^{n\times n}))$
with any $p>1$.

For any $t\in [0,T)$ and $\eta\in
L^2_{\cF_t}(\Om;\dbR^n)$, let us consider the following
forward-backward stochastic differential
equation:
\begin{equation}\label{5.14-eq12}
\left\{
\begin{array}{ll}\ds
d x^t(r)=\big(A+B\Th\big) x^tdr + \big(C +
D \Th \big)x^tdW(r) &\mbox{\rm in } [t,T],\\
dy^t(r) = -\big(A^{\top}  y^t + C^{\top}  z^t +
Qx^t \big) dr +  z^tdW(r) &\mbox{\rm in } [t,T],
\\ \ns\ds x^t(t)=\eta, \q
y^t(T)= G x^t(T).
\end{array}
\right.
\end{equation}
Clearly, equation \eqref{5.14-eq12} admits a
unique solution
$$
\big(x^t(\cd),y^t(\cd),z^t(\cd)\big)\in
L^2_\dbF(\Om;C([t,T];\dbR^n))\times
L^2_\dbF(\Om;C([t,T];\dbR^n))\times
L^2_\dbF(t,T;\dbR^n).$$
Also, consider the following forward-backward
stochastic differential equation:
\begin{equation}\label{5.14-eq12.1}
\left\{
\begin{array}{ll}\ds
d X^t(r)=\big(A+B\Th\big) X^tdr + \big(C +
D \Th \big)X^tdW(r) &\mbox{\rm in } [t,T],\\
dY^t(r) = -\big(A^{\top} Y^t + C^{\top} Z^t +
QX^t \big) dr + Z^tdW(r) &\mbox{\rm in } [t,T],
\\ \ns\ds X^t(t)=I_n, \q
Y^t(T)= G X^t(T)
\end{array}
\right.
\end{equation}
Likewise, equation \eqref{5.14-eq12.1} admits
a unique solution
$$
\big(X^t(\cd),Y^t(\cd),Z^t(\cd)\big)\in
L^2_\dbF(\Om;C([t,T];\dbR^{n\times n}))\times
L^2_\dbF(\Om;C([t,T];\dbR^{n\times n}))\times
L^2_\dbF(t,T;\dbR^{n\times n}).$$
It follows from \eqref{5.14-eq12} and
\eqref{5.14-eq12.1} that, for any $\eta\in L^2_{\cF_t}(\Om;\dbR^n)$,
\begin{equation}\label{7.19-eq6.1}
\begin{array}{ll}\ds
x^t(r) =  X^t(r)\eta, \q y^t(r) =  Y^t(r)\eta,\q& \forall\,r\in[t,T].
\\\ns\ds z^t(r) =  Z^t(r)\eta,&\ae r\in[t,T].
\end{array}
\end{equation}
By the uniqueness of the solution to
\eqref{5.7-eq4.1}, for any $\zeta\in \dbR^n$ and $t\in[0,T]$, we
have that
$$
X^t(r)X(t)\zeta = x^t(r;X(t)\zeta)=x(r;\zeta),\ \ \text{a.s.}
$$
thus,
$$
Y^t(t) X(t)\zeta = y^t(t;X(t)\zeta)=Y(t)\zeta. \ \ \text{a.s.}
$$
This implies that for all $t\in [0,T],$
\begin{equation}\label{5.14-eq13}
Y^t(t) =Y(t)\wt X(t)^{\top}= P(t).\q \dbP\hb{-a.s.}
\end{equation}
Let $\eta,\xi\in L^2_{\cF_t}(\Om;\dbR^n)$. Since
$Y^t(r) \eta=y^t(r;\eta)$ and $X^t(r)
\xi=x^t(r;\xi)$, applying It\^o's formula to
$\langle x^t(\cd),y^t(\cd) \rangle_{\dbR^n}$, we
get that
\begin{equation}\label{1.20-eq21}
\begin{array}{ll}\ds
\mE\langle \xi, P(t) \eta \rangle_{\dbR^n}
\3n&\ds= \mE\langle GX^t(T)\eta,
X^t(T)\xi \rangle_{\dbR^n} + \mE\int_t^T\langle Q(r)X^t(r)\eta,X^t(r)\xi \rangle_{\dbR^n} dr\\
\ns& \ds\q - \mE\int_t^T\langle B\Th X^t(r)\xi,
Y^t(r)\eta \rangle_{\dbR^n} dr-
\mE\int_t^T\langle D\Th X^t(r)\xi, Z^t(r)\eta
\rangle_{\dbR^n} dr.
\end{array}
\end{equation}
This, together with Corollary \ref{5.7-prop1}, implies
that
$$
\begin{array}{ll}\ds
\mE\langle P(t)\eta, \xi
\rangle_{\dbR^n}\3n&\ds=\mE\langle GX^t(T)
\eta,X^t(T)\xi \rangle_{\dbR^n}+\mE\int_t^T
\big(\langle
Q(r)X^t(r)\eta,X^t(r)\xi \rangle_{\dbR^n}\\
\ns&\ds \q + \langle
R(r)\Th(r)X^t(r)\eta,\Th(r)X^t(r)\xi
\rangle_{\dbR^n} \big)dr.
\end{array}
$$
Therefore,
\begin{equation}\label{5.14-eq14}
\begin{array}{ll}\ds
\mE\langle P(t)\eta, \xi \rangle_{\dbR^n}\3n&\ds = \mE\Big\langle \xi, X^{t}(T)^\top
GX^t(T)\eta\\
\ns&\ds \q +\mE \int_t^T \big(X^{t}(r)^\top
Q(r)X^t(r)\eta + X^{t}(r)^\top \Th(r)^{\top}
R(r)\Th(r)X^t(r)\eta \big)dr
\Big\rangle_{\dbR^n}.
\end{array}
\end{equation}
This concludes that
\begin{equation}\label{5.14-eq11}
\begin{array}{ll}\ds
P(t)  = \mE\(X^{t}(T)^\top GX^t(T)
\\
\ns\ds \qq \qq+\mE\int_t^T \big(X^{t}(r)^\top
Q(r)X^t(r) + X^{t}(r)^\top \Th(r)^{\top}
R(r)\Th(r)X^t(r) \big)dr\;\Big|\;\cF_t\).
\end{array}
\end{equation}
By \eqref{5.14-eq11} and the symmetry of $G$, $Q(\cd)$ and $R(\cd)$, it is easy to conclude that, for any $t\in[0,T]$, $P(t)$ is
symmetric, $\dbP$-a.s.

Next, we prove that $\L(t,\om)=\L(t,\om)^{\top}
$ for a.e. $(t,\om)\in (0,T)\times \Om$.

Clearly, $(P^\top,\L^\top)$ satisfies that
\begin{equation}\label{8.25-eq1}
\left\{
\begin{array}{ll}\ds
dP^{\top} =-\big[P^{\top} A + A^{\top} P^{\top}
+ \L^{\top}  C + C^{\top} \L^{\top}  +
C^{\top} P^{\top}  C \\
\ns\ds \qq\qq +\Th^{\top} (P B + C^{\top} PD +
\L D)^{\top}   + Q \big]dt + \L^{\top} dW(t)
&\mbox{ in }[0,T],\cr
\ns\ds P(T)^{\top}=G.
\end{array}
\right.
\end{equation}
 According to
\eqref{6.7-eq8} and \eqref{8.25-eq1}, and noting that $P(\cd)$ is symmetric, we find
that for any $t\in [0,T]$,
\begin{equation}\label{7.3-eq4}
\begin{array}{ll}\ds
0\3n&\ds=-\int_0^t \big\{\big[\L C+ C^{\top} \L
+(P
B + C^{\top} PD + \L D) \Th \big]\\
\ns&\ds \qq\qq -\big[\L C + C^{\top}  \L +(P B +
C^{\top} PD + \L D) \Th \big]^{\top} \big\}
d\tau + \int_0^t(\L-\L^{\top} ) dW(\tau).
\end{array}
\end{equation}
By \eqref{7.3-eq4} and the uniqueness of the
decomposition of semimartingale, we conclude
that
\bel{20160228e8}
\L(t,\om)=\L(t,\om)^{\top} ,\q \ae (t,\om)\in
(0,T)\times\Om.
\ee

\medskip

{\bf Step 3}. In this step, we show that
$(P,\L)$ is a pair of stochastic processes satisfying \eqref{5.5-eq6}, \eqref{20160225e2}, \eqref{5.7-eq5}, \eqref{5.10}. Moreover, \eqref{Value} holds.

From \eqref{5.7-eq5.1}, it holds that
\begin{equation}\label{5.7-eq6.1}
B^{\top} P+D^{\top} \Pi +R\Th =0, \q \ae
(t,\om)\in[0,T]\times\Om.
\end{equation}
By \eqref{5.7-eq6.1}, we see that
\begin{equation}\label{5.7-eq9}
\ba{ll}\ds
0=B^\top P + D^\top\big[\L + P(C+D\Th)\big] +
R\Th = B^\top P + D^\top PC + D^\top\L +K\Th\\\ns\ds\q=L+K\Th.
\ea
\end{equation}
Thus, it follows from \eqref{5.7-eq9} that $\cR(K(\cd))\supset\cR(L(\cd))$ and
$$
K^{\dag}K\Th=-K^{\dag}L.
$$
By Lemma \ref{lm3}, $K^{\dagger}K$ is an orthogonal projector. Hence we have
$$
\int_0^T|K^{\dag}(r)L(r)|_{\dbR^{m\times n}}^2dr= \int_0^T|K^{\dag}(r)K(r)|_{\dbR^{n\times n}}^2|\Th(r)|_{\dbR^{m\times n}}^2dr\leq
\int_0^T|\Th(r)|_{\dbR^{m\times n}}^2dr, \ \as
$$
This leads to \eqref{5.7-eq5}. Moreover, we have \eqref{5.10}, i.e.,
$\Th(\cd)=-K(\cd)^{\dag}L + \big(I_m -
K(\cd)^{\dag}K(\cd)\big)\th$
for some $\th\in
L^\infty_\dbF(\Om;L^2(0,T;\dbR^{m\times n}))$.
Therefore, by \eqref{20160228e8}, \eqref{5.7-eq9} and Lemma \ref{lm3}, it follows that
\begin{equation}\label{6.8-eq15}
\begin{array}{ll}
\ds (PB + C^\top PD + \L D )\Th=L^\top\Th =-\Th^\top K\Th
\\
\ns\ds =-\Th^\top K\left[-K(\cd)^{\dag}L + \big(I_m -
K(\cd)^{\dag}K(\cd)\big)\th\right]\\
\ns\ds =\Th^\top KK^{\dag}L= -L^\top K^{\dag}L.
\end{array}
\end{equation}
Hence, by \eqref{6.7-eq8}, we conclude that $(P,\L)$ is a solution to
\eqref{5.5-eq6}.

To obtain \eqref{20160225e2}, we only need to show that
\bel{20160226e6}
K\geq 0,\q\ae (t,\om)\in[0,T]\times\Om.
\ee
For this purpose, from  \eqref{6.8-eq15}, we see that
\begin{equation}\label{6.8dd-eq15}
\Th^\top K\Th =L^\top K^{\dag}L.
\end{equation}
Due to \eqref{6.8dd-eq15} and \eqref{9.7-eq10}, for any $(s,\eta)\in[0,T)\times L_{\cF_s}^2(\Omega;\dbR^n)$,  by repeating the procedures in deriving \eqref{5.31-eq1} above, we show that
\begin{equation}\label{7.3-eq11}
\ba{ll}\ds
\cJ(s,\eta;u(\cd))\3n&\ds=\frac{1}{2}\dbE\Big(\big\langle
P(s)\eta,\eta\big\rangle_{\dbR^n}+\int_s^T\big\langle K(u-\Th x),u-\Th x\big\rangle_{\dbR^m}dr\Big)\\
\ns&\ds
=\cJ\big(s,\eta;\Th(\cd)\bar x(\cd)\big)
+\frac{1}{2}\dbE\int_s^T\big\langle K(u-\Th
x),u-\Th x\big\rangle_{\dbR^m} dr.
\ea
\end{equation}
Hence, by the optimality of the feedback operator $\Th(\cd)$, \eqref{Value} holds and
\begin{equation}\label{8.22-eq1}
0\leq
\dbE\int_s^T\big\langle K(u-\Th
x),u-\Th x\big\rangle_{\dbR^m} dr,\qq\forall\; u(\cd)\in L^2_\dbF(s,T;\dbR^m).
\end{equation}
For any $v(\cd)\in L^2_\dbF(s,T;\dbR^m)$, we may choose a control $u(\cd)\in L^2_\dbF(s,T;\dbR^m)$ (in \eqref{5.2-eq1}) in the ``feedback form" $u(\cd)=v(\cd)+\Th (\cd)x(\cd)$. Hence, by \eqref{8.22-eq1}, we obtain \eqref{20160226e6}.
This completes the proof of
the necessity in Theorem \ref{5.7-th1}.

%%%%%%%%%%%%%%%%%%%%%%%%%%%%%%%%%%%%%%%%%%%%%%%%

\section{Proofs of Corollary \ref{corollary-1} and Theorem
\ref{5.7-th1-2}}\label{proof-2}

%%%%%%%%%%%%%%%%%%%%%%%%%%%%%%%%%%%%%%%%%%%%%%%%

This section is addressed to proving Corollary \ref{corollary-1} and Theorem
\ref{5.7-th1-2}.

\ms
{\it Proof of Corollary \ref{corollary-1}}\,: Suppose that the equation \eqref{5.5-eq6} admits two pairs of solution
$$\big(P_i(\cd),\L_i(\cd)\big)\in
L^{\infty}_{\dbF}(\Om;C([0,T];\cS(\dbR^n)))\times
L^p_{\dbF}(\Om;L^2(0,T;\cS(\dbR^n)))$$
($i=1, 2$), so
that
$$\ba{ll}
\ns\ds\cR(K_i(t,\omega))\supset\cR(L_i(t,\omega)),\ \  K_i(t,\omega)\geq 0,\ \ \ae (t,\omega)\in [0,T]\times\Omega,\\
\ns\ds K_i(\cd)^{\dag}L_i(\cd)\in
L^\infty_\dbF(\Om;L^2(0,T;\dbR^{m\times n})),
\ea$$
 where $K_i\= R+D^{\top}  P_iD$ and $
L_i\= B^{\top}  P_i+D^{\top}  (P_iC+\L_i)$. Let
$$\Th_i(\cd)\=-K_i(\cd)^{\dag}L_i(\cd) + \big(I_m -
K_i(\cd)^{\dag}K_i(\cd)\big)\th_i$$
 for some $\th_i\in
L^\infty_\dbF(\Om;L^2(0,T;\dbR^{m\times n}))$. Then by the sufficiency in Theorem \ref{5.7-th1}, $\Th_1(\cd)$ and $\Th_2(\cd)$ are two optimal feedback operators and
\begin{equation}\label{uniqueness-issue-1-2}
\inf_{u\in
L^2_\dbF(s,T;\dbR^m)}\cJ(s,\eta;u)=\frac{1}{2}\,\dbE\langle
P_1(s)\eta,\eta\rangle_{\dbR^n}=\frac{1}{2}\,\dbE\langle
P_2(s)\eta,\eta\rangle_{\dbR^n}.
\end{equation}
By the arbitrariness of $s,\ \eta$, one has $P_1(\cd)=P_2(\cd)$. Similar to \eqref{20160228e8}, one can show that $\L_1(\cd)=\L_2(\cd)$.\endpf

%Using the same tricks as in (\ref{6.8-eq15}) we have $L_i^{\top}\Th_i=-L_i^{\top}K_i^{\dagger}L_i$. Since $(P_0,\L_0)$ solves the Riccati equation  \eqref{5.5-eq6}, we see that
%$(P_0(\cd),\L_0(\cd))$ solves the following
%$\dbR^{n\times n}$-valued backward stochastic
%differential equation:
%%
%$$
%\left\{
%\begin{array}{ll}\ds
%dP_0 =-\big[P_0A+ A^{\top}P_0 + \L_0 C+ C^{\top}\L_0+
%C^{\top} P_0 C \\
%\ns\ds \qq\q +(P_0 B + C^{\top}P_0D + \L_0 D) \Th + Q
%\big]dt + \L_0 dW(t) &\mbox{ in }[0,T],\cr
%%
%\ns\ds P_0(T)=G,
%\end{array}
%\right.
%$$
%%
%which is the same as \eqref{6.7-eq8}. By the uniqueness result in Lemma \ref{lm2.1}, we conclude that $(P_0(\cd),\L_0(\cd))=(P(\cd),\L(\cd))$.

\ms

\ms

{\it Proof of  Theorem \ref{5.7-th1-2}}\,: {\bf
The ``if" part}. By the necessity in  Theorem
\ref{5.7-th1}, it remains to show the uniqueness
of optimal feedback operators. Suppose there
exists another optimal feedback operator $
\wt\Th(\cd)$. By the necessity in Theorem
\ref{5.7-th1}, the Riccati equation
\eqref{5.5-eq6} admits a unique solution
$$\big(\wt P(\cd),\wt \L(\cd)\big)\in
L^{\infty}_{\dbF}(\Om;C([0,T];\cS(\dbR^n)))
\times L^p_{\dbF}(\Om;L^2(0,T;\cS(\dbR^n)))$$
so
that
$$\ba{ll}
\ns\ds \cR(\wt K(t,\omega))\supset\cR(\wt
L(t,\omega)),\ \ \wt K(t,\omega)\geq 0, \ \ \ae
(t,\omega)\in
[0,T]\times\Omega,\\
 \ns\ds \wt K(\cd)^{\dag}\wt L(\cd)\in
L^\infty_\dbF(\Om;L^2(0,T;\dbR^{m\times n})),\ \
\wt \Th(\cd)=-\wt K(\cd)^{\dag}\wt L(\cd) +
\big(I_m - \wt K(\cd)^{\dag}\wt K(\cd)\big)\wt
\th
\ea$$
 for some $\wt \th \in
L^\infty_\dbF(\Om;$ $L^2(0,T;\dbR^{m\times
n}))$, where $\wt K \= R+D^{\top}\wt P D$ and $
\wt L \= B^{\top}  \wt P +D^{\top}  (\wt P C+\wt
\L)$. Moreover,
\begin{equation}\label{Value2016}
\inf_{u\in L^2_\dbF(s,T;\dbR^m)}\cJ(s,\eta;u)=
\frac{1}{2}\,\dbE\langle \wt P
(s)\eta,\eta\rangle_{\dbR^n}.
\end{equation}
Since \eqref{Value} and \eqref{Value2016} hold
for  any $(s,\eta)\in [0,T)\times
L^2_{\cF_s}(\Om;\dbR^n)$, it follows that
$P(\cd)=\wt P(\cd)$. Similar to
\eqref{20160228e8}, one can show that
$\L(\cd)=\wt \L(\cd)$. Hence, $K=\wt K$ and
$L=\wt L$. Since $K(t,\omega)> 0$ for $\ae
(t,\omega)\in
 [0,T]\times\Omega$, one has $\Th(\cd)=\wt \Th(\cd)$.

\ms

{\bf The ``only if" part}. We only need to prove that the uniqueness and existence of optimal feedback operators implies $K(\cd)>0$ a.e. For any $\tilde\th\in L^{\infty}_{\dbF}(\Omega;L^2(0,T;\dbR^{m\times n}))$, we construct another stochastic process $\wt\Th\in L^{\infty}_{\dbF}(\Omega;L^2(0,T;\dbR^{m\times n}))$ as follows
$$
\wt\Th\triangleq -K^{\dagger}L+(I_m-K^{\dagger}K)(\th+\tilde\th).
$$
Repeating the argument in the proof of sufficiency in Theorem \ref{5.7-th1}, one can show that $\wt\Th(\cd)$ is an optimal feedback operator. By the uniqueness of optimal feedback operators, we deduce that $\Th(\cd)=\wt\Th(\cd)$, and therefore $(I_m-K^{\dagger}K)\th'=0$. The arbitrariness of $\th'$ indicates that $K^{\dagger}K=I_m$. As a result, $K^{\dagger}=K^{-1}$, and hence $K(\cd)>0$ a.e.\endpf
%
%By the sufficiency in Theorem \ref{5.7-th1}, we have $K(\cd)\geq0$ a.e.  From

%%%%%%%%%%%%%%%%%%%%%%%%%%%%%%%%%%%%%%%%%%%%%%%%%%%%

\section{Two illustrating
examples}\label{example}

%%%%%%%%%%%%%%%%%%%%%%%%%%%%%%%%%%%%%%%%%%%%%%%%%%%%

We have discussed the relationship between the
existence of feedback operator  and the
well-posedness of the Riccati equation
\eqref{5.5-eq6}. In this section, we give two
examples which are  inspired by
\cite{Wang-2016}. In the  first example, we show that there is a
feedback operator $\Th(\cd)\in
L^\infty_\dbF(\Om;L^2(0,T;\dbR^{m\times n}))$; while in the second one, it is shown that the desired feedback operator does
not exist.

\begin{example}
Applying It\^o's formula to $\sin W(\cd)$, we obtain that
 \bel{201611191}
 \sin W(T)-\sin W(t)=\int_t^T\cos W(s)dW(s)-\frac{1}{2}\int_0^T\sin W(s)ds.
 \ee
Write
 \bel{201611192}
 \ba{ll}
  \ds \xi\=2+\frac{T}{2}+\sin W(T)+\frac{1}{2}\int_0^T\sin W(s)ds,\\[3mm]\ds y(t)\=2+\frac{T}{2}+\sin W(t)+\frac{1}{2}\int_0^t\sin W(s)ds,\quad Y(t)\=\cos W(t),\qq t\in[0,T].
  \ea
 \ee
From \eqref{201611191}, it is clear that $1\le y(\cd)\le 3+T$, and $(y(\cd),Y(\cd))$ satisfies
$$ y(t)=\xi-\int_t^TY(s)dW(s),\qq t\in[0,T].
$$

Consider an SLQ problem with the following data (Note that, by \eqref{201611192}, $1\le \xi\le 3+T$):
\bel{201611193}
m=n=1,\ \ A=B=C=Q=S=0,\ \ D=1,\ \ R=\frac{1}{2(3+T)},\ \ G=\xi^{-1}-R>0.
\ee
The corresponding Riccati equation is
\begin{equation}\label{Riccati-equation-example}
\left\{\ba{ll}
dP(s)=(R+P(s))^{-1}\L^2(s)ds+\L(s)dW(s),\q s\in [0,T],\\[3mm]
P(T)=G.
\ea\right.
\end{equation}
By It\^o's formula, one can prove that $\big(P(\cd),\L(\cd)\big)=\big(y(\cd)^{-1}-R,-y(\cd)^{-2}Y(\cd)\big)$
is the unique solution to
\eqref{Riccati-equation-example}. According to Theorem
\ref{5.7-th1},
$$\Th(\cd)\triangleq
-y(\cd)^{-1}Y(\cd)\in
L^{\infty}_{\dbF}(\Omega;L^2(0,T;\dbR))$$
is an
optimal feedback operator.

\end{example}

Next, we give an negative example to show the
nonexistence of the optimal feedback operator.
\begin{example}\label{counterexample-1}
Define one-dimensional stochastic processes $M(\cd),\ \zeta(\cd)$ and
stopping time $\t$ as follows:
\begin{equation}\label{11.12-eq1}
\begin{cases}
\ds M(t)\triangleq \int_0^t\frac{1}{\sqrt{T-s}}dW(s),\qq t\in[0,T),\\
\ns\ds \t\triangleq \inf\big\{t\in[0,T),\ |M(t)|>1\big\}\wedge T,\\
\ns\ds \zeta(t)\triangleq
\frac{\pi}{2\sqrt{2}\sqrt{T-t}}\chi_{[0,\t]}(t),\qq
t\in[0,T).
\end{cases}
\end{equation}
It was shown in \cite[Lemma A.1]{FreidosReis}
that
\bel{boundedness-example}\ba{ll}
\ns\ds
\Big|\int_0^T\zeta(s)dW(s)\Big|=\frac{\pi}{2\sqrt{2}}\Big|\int_0^\t\frac{1}{
\sqrt{T-t}}
dW(t)\Big|=\frac{\pi}{2\sqrt{2}}\big|M(\t)\big|\leq
\frac{\pi}{2\sqrt{2}},
\ea\ee
and
\bel{aim-1}\ba{ll}
\ns\ds
\dbE\Big[\exp\Big(\int_0^T|\zeta(t)|^2dt\Big)\Big]=\infty.
\ea\ee

Consider the following backward stochastic differential equation:
$$
Y(t)=\int_0^T\zeta(s)dW(s)+\frac{\pi}{2\sqrt{2}}+1-\int_t^TZ(s)dW(s),\qq
t\in[0,T].
$$
This equation admits a unique solution $(Y,Z)$ as follows
$$
Y(t)=\int_0^t\zeta(s)dW(s)+\frac{\pi}{2\sqrt{2}}+1,\quad
 Z(t)=\zeta(t),\qquad t\in[0,T].
$$
From \eqref{11.12-eq1}--\eqref{aim-1}, it is easy to see that
\begin{equation}\label{11.12-eq2}
\begin{cases}
1\leq Y(\cd)\leq \frac{\pi}{\sqrt{2}}+1, \\
Z(\cd) \notin
L^{\infty}_{\dbF}(\Omega;L^2(0,T;\dbR)).
\end{cases}
\end{equation}

Consider an SLQ problem with the following data:
\bel{201611195}
m=n=1, \ \ A=B=C=Q=S=0,\ \ D=1,\ \ R=\frac{1}{4}>0,\ \
G=Y(T)^{-1}-\frac{1}{4}>0.
\ee
For this problem, the corresponding Riccati equation reads
\begin{equation}\label{Riccati-equation-example00}
\left\{\ba{ll}
dP(s)=(R+P(s))^{-1}\L^2(s)ds+\L(s)dW(s),\q s\in [0,T],\\[3mm]
P(T)=G,
\ea\right.
\end{equation}
and $\Th(\cd)=-(R+P(\cd))^{-1}\L(\cd)$.

Put
$$
\wt P(\cd)\triangleq P(\cd)+R,\ \ \wt
\L\triangleq \L.
$$
It follows from \eqref{Riccati-equation-example00}
that
\begin{equation}\label{Riccati-equation-example-wt}\left\{\ba{ll}
d\wt P(s)=\wt
P(s)^{-1}\wt\L^2(s)ds+\wt\L(s)dW(s),\q s\in [0,T],\\[3mm] \wt
P(T)=Y(T)^{-1}.
\ea\right.
\end{equation}
Applying It\^{o}'s formula to $Y(\cd)^{-1}$, we
deduce that $(\wt P(\cd),\wt \L(\cd))=(Y(\cd)^{-1},-Y(\cd)^{-2}Z(\cd))$ is
the unique solution to
\eqref{Riccati-equation-example-wt}. As a
result,
$$ (P(\cd),\L(\cd))\triangleq ( Y(\cd)^{-1}-R,-Y(\cd)^{-2}Z(\cd))
$$
is the unique solution to the Riccati equation
(\ref{Riccati-equation-example00}). Moreover,
$\Th(\cd)=-Y(\cd)^{-1}Z(\cd)$. By
\eqref{11.12-eq2}, we see that $\Th(\cd)$
does not belong to
$L^{\infty}_{\dbF}(\Omega;L^2(0,T;\dbR))$, either. Hence, it is not a ``qualified" feedback operator.
\endpf
\end{example}

\br
Clearly, the form of \eqref{Riccati-equation-example00} is the same as that of \eqref{Riccati-equation-example} but their endpoint values at $T$ are different. For the endpoint value $G$ given in \eqref{201611192}, the corresponding $\L(\cd) \in
L^{\infty}_{\dbF}(\Omega;L^2(0,T;\dbR))$. However, for the endpoint value $G$ given in \eqref{201611195}, the resulting $\L(\cd) \notin
L^{\infty}_{\dbF}(\Omega;L^2(0,T;\dbR))$.

Generally speaking, it would be quite interesting to find some suitable conditions to guarantee that the equation
\eqref{5.5-eq6} admits a unique solution
$(P(\cd),\L(\cd))\in L^\infty_\dbF(0,T;\cS(\dbR^n))\times
L^\infty_\dbF(\Om;L^2(0,T;\cS(\dbR^n)))$ but this is an unsolved problem.
\er

\br
Example \ref{counterexample-1} also shows that, a solvable Problem (SLQ) does not need to have feedback controls. This is a significant difference between SLQs and their deterministic counterparts. Indeed, it is well-known that one can always find the desired feedback control through the corresponding Riccati equation whenever a deterministic LQ is solvable.
\er

%%%%%%%%%%%%%%%%%%%%%%%%%%%%%%%%%%%%%%%%%%%%%%%%%%%%%%%

\section*{Acknowledgement}

%%%%%%%%%%%%%%%%%%%%%%%%%%%%%%%%%%%%%%%%%%%%%%%%%%%%%%%

This work is supported by the NSF of China under
grants 11471231, 11231007, 11301298 and 11401404, the PCSIRT under grant IRT$\!\_$15R53 and the Chang Jiang
Scholars Program from Chinese Education
Ministry, and the Fundamental Research Funds for
the Central Universities in China under grant
2015SCU04A02. The authors gratefully acknowledge Professor Jiongmin Yong for helpful discussions.


\begin{thebibliography}{99}



\bibitem{AMZ} M.~Ait Rami, J.~B.~Moore and X.~Zhou. \it
Indefinite stochastic linear quadratic control
and generalized differential Riccati equation.
\sl SIAM J. Control Optim. \rm {\bf 40} (2001),
1296--1311.

\bibitem{Athens} M.~Athans. \it The role and use of the
stochastic linear-quadratic-Gaussian problem in
control system design. \sl IEEE Trans. Automat.
Control. \rm {\bf 16} (1971), 529--552.

\bibitem{Ben-Israel} A.~Ben-Israel and T.~N.~E.~Greville.
\sl Generalized Inverses: Theory and
Applications. \rm Pure and Applied Mathematics.
Wiley-Interscience [John Wiley \& Sons], New
York-London-Sydney, 1974.

\bibitem{Bensoussan1} A.~Bensoussan. \it Lectures on stochastic control. \sl Nonlinear Filtering and Stochastic Control, \rm 1--62. Lecture Notes in Math., vol. {\bf 972}. Springer-Verlag, Berlin, 1981.

\bibitem{Bismut1} J.-M.~Bismut. \it Linear quadratic optimal
stochastic control with random coefficients. \sl
SIAM J. Control Optim. \rm{\bf 14} (1976),
419--444.

\bibitem{Bismut2} J.-M.~Bismut. \it  Contr\^ole des syst\`emes
lin\'eaires quadratiques: applications de
l'int\'egrale stochastique. \sl S\'eminaire de
Probabilit\'es XII, Universit\'e de Strasbourg
1976/77, \rm 180--264. Lecture Notes in Math.,
\rm  vol. {\bf 649}, Springer-Verlag, Berlin,
1978.


\bibitem{BDHPS}Ph. Briand, B. Delyon, Y. Hu, E. Pardoux and L. Stoica. \it $L^p$ solutions of backward stochastic differential
equations. \sl Stochastic Process. Appl. \rm {\bf 108} (2003), 109--129.

\bibitem{CLZ1} S.~Chen, X.~Li and X.~Zhou. \it
Stochastic linear quadratic regulators with
indefinite control weight costs. \sl SIAM J.
Control Optim. \rm{\bf 36} (1998), 1685--1702.

\bibitem{Davis} M.~H.~A.~Davis. \sl Linear Estimation and Stochastic
Control. \rm Chapman and Hall Mathematics
Series. Chapman and Hall, London; Halsted Press
[John Wiley $\&$ Sons], New York, 1977.

\bibitem{Delbaen-Tang} F. Delbaen and S. Tang. \it Harmonic analysis of stochastic equations and backward stochastic differential equations. \sl Probab. Theory Relat. Fields.  \rm{\bf 146} (2010), 291--336.

\bibitem{FreidosReis} C. Frei and G. dos Reis. \it A financial market with interacting investors: does an equilibrium
exist? \sl Math. Finan. Econ. \rm{\bf 4 } (2011), 161--182.

\bibitem{Kalman} R.~E.~Kalman. \it Contributions to the theory of
optimal control. \sl Bol. Soc. Mat. Mexicana.
\rm {\bf 5} (1960), 102--119.

\bibitem{LZ0}  Q. L\"u and X.~Zhang. \it Well-posedness of backward stochastic differential equations with
general filtration. \sl J.  Differential Equations. \rm {\bf 254}
(2013), 3200--3227.

\bibitem{LZ1} Q.~L\"u and X.~Zhang. \sl General Pontryagin-Type
Stochastic Maximum Principle and Backward
Stochastic Evolution Equations in Infinite
Dimensions. \rm Springer Briefs in Mathematics.
Springer, Cham, 2014.

\bibitem{LZ} Q.~L\"u and X.~Zhang. \it Transposition
method for backward stochastic evolution
equations revisited, and its application. \sl
Math. Control Relat. Fields. \rm{\bf 5} (2015),
529--555.

\bibitem{LZ3}  Q. L\"u and X.~Zhang. \it Optimal feedback for stochastic linear quadratic control and backward stochastic Riccati
equations in infinite dimensions. \rm Preprint.

\bibitem{PP} E.~Pardoux and S.~Peng. \it Adapted solution of backward
stochastic equation. \sl Systems Control Lett. \rm {\bf 14} (1990),
55--61.

\bibitem{Peng} S.~Peng. \it Stochastic Hamilton-Jacobi-Bellman
equations. \sl SIAM J. Control Optim. \rm{\bf
30} (1992), 284--304.

\bibitem{Huyen} H.~Pham. \it Linear quadratic optimal control of conditional McKean-Vlasov equation with random coefficients and applications. \rm arXiv:1604.06609v1.

\bibitem{Protter} P.~E.~Protter. \sl Stochastic Integration and
Differential Equations. \rm Stochastic Modelling
and Applied Probability, vol. {\bf 21}. Springer-Verlag,
Berlin, 2005.

\bibitem{Reid} W.~T.~Reid. \it A matrix differential equation
of Riccati type. \sl Amer. J. Math. \rm{\bf 68}
(1946), 237--246.

\bibitem{SY} J.~Sun and J.~Yong. \it Linear quadratic
stochastic differential games: open-loop and
closed-loop saddle points. \sl SIAM J. Control
Optim. \rm{\bf 52} (2014), 4082--4121.

\bibitem{Tang1} S.~Tang. \it  General linear quadratic optimal
stochastic control problems with random
coefficients: linear stochastic Hamilton systems
and backward stochastic Riccati equations. \sl
SIAM J. Control Optim. \rm{\bf 42} (2003), 53--75.

\bibitem{Tang2} S.~Tang. \it Dynamic programming for general
linear quadratic optimal stochastic control with
random coefficients. \sl SIAM J. Control Optim.
\rm {\bf 53} (2015), 1082--1106.

\bibitem{Wonham2} W.~M.~Wonham. \it On a matrix Riccati equation of stochastic control. \sl SIAM J.
Control. \rm {\bf 6} (1968), 681--697.

\bibitem{Wonham1} W.~M.~Wonham. \sl Linear Multivariable Control,
a Geometric Approach. \rm Applications of
Mathematics, vol. {\bf 10}. Springer-Verlag, New
York, 1985.

\bibitem{Wang-2016} T. Wang. \it New optimality conditions in linear quadratic problems with random coefficients and applications. \rm In submission.


\bibitem{YL} J.~Yong  and   H.~Lou. \sl A Concise Course on  Optimal Control Theory. \rm Higher Education Press, Beijing, 2006. (In Chinese)

\bibitem{YZ}J. Yong and X.Y. Zhou. \sl Stochastic Controls: Hamiltonian Systems and HJB
Equations. \rm Springer-Verlag, New York, Berlin, 2000.


\end{thebibliography}
\end{document}